\documentclass[12pt]{article}
\usepackage[english]{babel}
\usepackage{amssymb,amsthm}
\usepackage{amsmath}
\usepackage{inputenc}
\usepackage{graphicx}
\newtheorem{thm}{Theorem}[section]

\newtheorem{pr}{Proposition}[section]
\newtheorem{Def}{Definition}[section]

\begin{document}

\begin{center}
{\large\bf Classifications of some classes of Zinbiel algebras.\\[2mm]}
\end{center}

\begin{center}
{\bf J.Q. Adashev \footnote{ This work is supported by ICTP}\\
{\small\it Department of Mathematics, National University of
Uzbekistan, Tashkent, Uzbekistan \\
e-mails: adashevjq@mail.ru}}
\end{center}

\begin{center}{ \bf B.A. Omirov \footnote{ and DFG 436 USB
 113/10/0-1 project (Germany).}\\
{\small\it Institute of Mathematics, Uzbek Academy of Sciences, Tashkent, Uzbekistan \\
e-mail: omirovb@mail.ru}}
\end{center}

\begin{center}
{\bf  A.Kh. Khudoyberdiyev\\
{\small\it Department of Mathematics, National University of
Uzbekistan, Tashkent, Uzbekistan \\
e-mail:  abror\_alg\_85@mail.ru}}
\end{center}

\begin{abstract} In this work nul-filiform and filiform Zinbiel algebras
are described up to isomorphism. Moreover, the classification of
complex Zinbiel algebras is extended from dimensions $\leq 3$ up
to the dimension $4.$


{\bf MR(2000) Subject Classification:} {\rm 17A32.}
\end{abstract}

\section{Introduction}

One of the important object of the modern theory of
non-associative algebras is Lie algebras. Active investigations in
the theory of Lie algebras lead to the appearing of some
generalizations of these algebras such as Mal'cev algebras, Lie
superalgebras, binary Lie algebras, Leibniz algebras and others.

In the present work, we consider algebras which are dual to
Leibniz algebras. Recall that Leibniz algebras were introduced in
\cite{Lod2} in the 90-th years of the last century. They are
defined by the following identity
$$[x,[y,z]]=[[x,y],z] - [[x,z],y].$$

J.-L. Loday in \cite{Lod1} studied categorical properties of
Leibniz algebras and considered in this connection a new object
--- Zinbiel algebra (read Leibniz in reverse order). Since the
category of Zinbiel algebras is Koszul dual to the category of
Leibniz algebras, sometimes they are also called dual Leibniz
algebras \cite{Lod3}.

In works \cite{Dzhu1}--\cite{Dzhu2} some interesting properties of
Zinbiel algebras were obtained. In particular, the nilpotency of
an arbitrary complex finite dimensional Zinbiel algebra was proved
in \cite{Dzhu2}.

For the examples of Zinbiel algebras we refer to works
\cite{Dzhu2}, \cite{Lod1} and \cite{Lod3}.

Since description of all finite dimensional complex Zinbiel
algebras (which are nilpotent) is a boundless problem, it is
natural to add certain restrictions for their investigation. One
of such restriction is restriction on the nilindex. Note works
\cite{Ayu} and \cite{Ver}, where structural description of
nilpotent Lie algebras and Leibniz algebras is given.

In study of any class of algebras, it is important to describe up
to isomorphism even algebras of lower dimensions because such
description gives examples for to establish or reject certain
conjectures. In this way in \cite{Dzhu2}, the classification of
complex Zinbiel algebras of dimensions $\leq 3$ is given. Applying
some general results obtained for finite dimensional Zinbiel
algebras we extended classification of complex Zinbiel algebras up
to dimension $4.$

\section{Classification of complex nul-filiform Zinbiel algebras}

\begin{Def}
An algebra $A$ over a field $F$ is called {\sl Zinbiel algebra} if
for any $x, y, z \in A$  the following identity
\begin{equation}\label{e1}
(x\circ y)\circ z=x\circ(y\circ z)+x\circ (z\circ y)
\end{equation}
holds. \end{Def}

For a given Zinbiel algebra $A$ we define the following sequence:
$$A^1=A, A^{k+1}=A\circ A^k, k \geq 1.$$

\begin{Def} An Zinbiel algebra $A$ is called {\sl nilpotent} if
there exists $s \in N$ such that $A^s=0$. The minimal number $s$
satisfying this property is called {\sl index of nilpotency} or
{\sl nilindex} of the algebra $A$.
 \end{Def}

It is not difficult to see that the index of nilpotency of an
arbitrary $n$-dimensional nilpotent algebra does not exceed the
number $n+1.$

\begin{Def} An $n-$dimensional Zinbiel algebra $A$ is called
{\sl nul-filiform} if $dim A^i=(n+1)-i,$ $1 \leq i \leq n+1$.
\end{Def}

It is evident, the last definition is equivalent to the fact that
algebra $A$ has maximal index of nilpotency.

\begin{thm}\label{t1} An arbitrary $n$-dimensional nul-filiform Zinbiel
algebra is isomorphic to the algebra:
\begin{equation}\label{e2}
e_i \circ e_j=C_{i+j-1}^j e_{i+j}, \ \mbox{ for } \ 2\leq i+j \leq
n
\end{equation}
where omitted products are equal to zero and $\{e_1, e_2, \ldots,
e_n\}$ is a basis of the algebra.
\end{thm}
\begin{proof}
Let $A$ be an $n-$dimensional nul-filiform Zinbiel algebra and let
$\{x_1, x_2,  \\ \ldots , x_n\}$ is a basis of the algebra $A$
such that $x_1 \in A^1 \setminus A^2$, $x_2\in A^2\setminus A^3$,
$\dots,$  $x_n \in A^n$. Since $x_2 \in A^2\setminus A^3$, we have
for some elements $b_{2,p}, c_{2,p}$ of algebra $A$
$$
x_2=\sum{b_{2,p}\circ c_{2,p}}= \sum{\alpha_{i,j} x_i\circ
x_j}=\alpha_{1,1} x_1\circ x_1+(*)
$$
where $(*) \in A^3$, i.e. $x_2=\alpha_{1,1}x_1 \circ x_1+(*)$.
Note that $\alpha_{1,1}x_1 \circ x_1 \neq0$. Indeed, in the
opposite case $x_2 \in A^3$.

Similar, for  $x_3 \in A^3\setminus A^4$ we have
$$
x_3= \sum a_{3,p}\circ (b_{3,p} \circ c_{3,p})= \sum
\alpha_{i,j,k} x_i\circ (x_j\circ x_k) =\alpha_{1,1,1} x_1\circ
(x_1\circ x_1)+(**)
$$
where $(**) \in A^4$ and $\alpha_{1,1,1}x_1\circ(x_1\circ x_1)\neq
0$ (otherwise $x_3\in A^4$), i.e. $x_3=\alpha_{1,1,1} x_1\circ
(x_1\circ x_1)+(**)$. Continuing this procedure we obtain that
elements $e_1:=x_1,$ $e_2:= x_1\circ x_1,$ $e_3:= x_1\circ
(x_1\circ x_1),$ $\dots,$ $e_n:= (x_1\circ \ldots \circ (x_1\circ
(x_1\circ x_1)))$ are different from zero. It is not difficult to
check the linear independence of these elements. Hence we can
choose as a basis of algebra $A$ the elements $\{e_1, e_2, \ldots
, e_n\}$. We have by construction
\begin{equation}\label{e3}
e_1\circ e_i=e_{i+1}\  \mbox{ for } \  1 \leq i \leq n-1.
\end{equation}

We shall prove the equality (\ref{e2}) by induction on $j$ for any
$i$.

Using identities (1), (3) we can prove by induction the equality:
$$
e_i\circ e_1=ie_{i+1} \  \mbox{ for } \   1 \leq i \leq
n-1,
$$
i.e. equality (\ref{e2}) is true for  $j=1$ and any $i$.

Suppose that the equality is true for all $j \leq k-1$ and any
$i$.

Let us prove the equality (\ref{e2}) for $j=k$ and any $i$. Using
the inductive hypothesis and the following chain of equalities:
$$
e_i\circ e_k= e_i\circ (e_1\circ e_{k-1})=(e_i\circ e_1)\circ
e_{k-1}- e_i\circ (e_{k-1}\circ e_1)=
$$ $$
=ie_{i+1}\circ
e_{k-1}-(k-1)e_i\circ e_k=i C_{i+k-1}^{k-1} e_{i+k} -(k-1)e_i\circ
e_k
$$
we obtain $ ke_i\circ e_k=iC_{i+k-1}^{k-1} e_{i+k}, $ i.e.
$$
e_i\circ e_k=\frac {i} {k} C_{i+k-1}^{k-1} e_{i+k}= \frac {i} {k}
\frac{(i+k-1)!} {(k-1)! i!} e_{i+k}= \frac{(i+k-1)!} {k! (i-1)!}
e_{i+k}= C_{i+k-1}^ke_{i+k}.
$$ \end{proof}

We denote the algebra from theorem \ref{t1} as $NF_n$.

It is not difficult to see that $n$-dimensional Zinbiel algebra is
one generated if and only if it is isomorphic to the algebra
$NF_n$.

\section{Classification of complex filiform Zinbiel algebras}

\begin{Def} An $n-$dimensional Zinbiel algebra $A$ is said to be {\sl filiform} if
$dim A^i=n-i,$ $2 \leq i \leq n$.
\end{Def}

Consider the natural gradation of an filiform Zinbiel algebra $À$
taking $A_i=A^i/A^{i+1},$ $1 \le i \le n-1$. It is obvious that
$dim A_1=2$ and $dim A_i=1,$ $2 \leq i \leq n-1,$ and
$$
A\cong A_1\oplus A_2 \oplus \ldots \oplus A_{n-1} \
\hbox{ where} \ A_i \circ A_j\subseteq A_{i+j}.
$$

In the following theorem, the classification of complex naturally
graded filiform Zinbiel algebras is represented.

\begin{thm}\label{t2} An arbitrary $n$-dimensional $(n\geq 5)$ naturally graded
complex filiform Zinbiel algebra is isomorphic to the following
algebra:
\begin{equation}\label{e4}
F_n: e_i\circ e_j= C_{i+j-1}^j e_{i+j}, \quad  2\leq i+j \leq n-1
\end{equation}
where omitted products are equal to zero and $\{e_1, e_2, \ldots,
e_n\}$ is a basis of the algebra.
\end{thm}
\begin{proof} Let $A$ be an Zinbiel algebra satisfying conditions of the
theorem. Similar to the work \cite{Ver} we take a basis $\{e_1,
e_2, \dots, e_n\}$ of the algebra $A$ such that $A_1=\langle e_1,
e_n\rangle,$ $A_2=\langle e_2\rangle,$ $A_3=\langle e_3\rangle,$
$\dots,$ $A_{n-1}=\langle e_{n-1}\rangle,$ and $e_1 \circ
e_i=e_{i+1}$ at $2\leq i\leq n-2$.

Introduce the following notations
$$
e_1\circ e_1=\alpha_1e_2,\quad e_1\circ e_n=\alpha_2e_2,\quad
e_n\circ e_1=\alpha_3e_2,\quad e_n\circ e_n=\alpha_4e_2.
$$

Now consider the possible cases.

\textbf{Case 1.} Let $(\alpha_1, \alpha_4)\neq (0, 0)$. Then
without any loss of generality we can suppose $\alpha_1\neq 0.$
Change the basis as follows $e'_2=\alpha_1e_2,$ $e'_3=\alpha_1
e_3,$ $\dots,$ $e'_{n-1}=\alpha_1e_{n-1}$ we can suppose
$\alpha_1=1$.

In this case, the space spanned on vectors $\{e_1, e_2, \dots,
e_{n-1}\}$ forms an nul-filiform Zinbiel algebra of the dimension
$n-1$. From the proof of theorem \ref{t1} we can conclude
$e_i\circ e_j= C_{i+j-1}^j e_{i+j},$ $2\leq i+j \leq n-1.$

Consider identity (\ref{e1}) in the following multiplications:
$$
(e_1\circ e_n)\circ e_1=e_1\circ (e_n\circ e_1)+ e_1\circ
(e_1\circ e_n) \Rightarrow 2\alpha_2e_3=\alpha_3e_3+\alpha_2e_3,
$$
i.e. $\alpha_2=\alpha_3;$
$$
(e_1\circ e_1)\circ e_n=e_1\circ (e_1\circ e_n)+ e_1\circ(e_n\circ
e_1) \Rightarrow  e_2\circ e_n=\alpha_3e_3+\alpha_2e_3,
$$
hence, $e_2\circ e_n=2\alpha_2e_3;$
$$
(e_1\circ e_n)\circ e_n=2e_1\circ (e_n\circ e_n) \Rightarrow 2
\alpha_2^2e_3=2\alpha_4e_3,
$$
consequently, $\alpha_2^2=\alpha_4;$
$$(e_n\circ e_1)\circ e_1=2e_n\circ (e_1\circ e_1) \Rightarrow
2\alpha_2e_3=2e_n\circ e_2,
$$
i.e. we have $e_n\circ e_2=\alpha_2e_3.$

Take the change of the basic elements by the following way:
$$
e'_n=-\alpha_2e_1+e_n,\  e'_i =e_i \  \mbox{\rm for} \ 1 \leq i
\leq n-1.
$$

It is not difficult to see that $e'_1\circ e'_n=e'_n\circ
e'_1=e'_n\circ e'_n=e'_2\circ e'_n=e'_n\circ e'_2=0,$ moreover
other products are not changed, i.e. we can suppose
$\alpha_2=\alpha_3=\alpha_4=0$.

Using identity (\ref{e1}) and the method of mathematical induction
one can easy to prove
\begin{equation}\label{e5}
e_n\circ e_i=0\  \mbox{ for } \ 1 \leq i \leq n-1.
\end{equation}

The equality
$$
e_i\circ e_n=0,\ \mbox{ for } \ 1 \leq i \leq n-1
$$
can be proved by induction using (\ref{e1}) and (\ref{e5}).

\textbf{Case 2.} Let $(\alpha_1, \alpha_4)=(0, 0)$. Then
$(\alpha_2, \alpha_3)\ne (0, 0)$. In the case of $\alpha_2 \ne
-\alpha_3$ we obtain case 1 taking $e'_1=e_1+e_n$. Now consider
the case $\alpha_2=-\alpha_3\neq 0$. Taking the following change
of basis
$$
e'_1 =e_1, \ \  e'_n=e_n, \ \ e'_i=\alpha_2e_i \ \mbox{ for } \ 2
\leq i \leq n-1
$$
we can suppose $\alpha_2=1$.

Consider the following products
$$
(e_1\circ e_n)\circ e_1=e_1\circ (e_n \circ e_1)+ e_1\circ
(e_1\circ e_n) \Rightarrow e_2\circ e_1=0;
$$
$$
0=(e_1\circ e_1)\circ e_2=e_1\circ (e_1\circ e_2)+ e_1\circ
(e_2\circ e_1)=e_1\circ e_3=e_4.
$$
Thus, we obtain the contradiction with existence of an algebra in
this case.
\end{proof}

The following proposition allows to extract a "convenient"\ basis
in an arbitrary complex filiform Zinbiel algebra. Such basis is
often called in literature adapted \cite{Ver}.

\begin{pr}\label{p1} There exists basis $\{e_1, e_2, \ldots, e_n\}$
in an arbitrary $n$-dimen\-sional $(n\geq 5)$ complex filiform
Zinbiel algebra such that the multiplication of the algebra has
the following form:
\begin{equation}\label{e6}
\begin{array}{c} e_i\circ e_j=C_{i+j-1}^j e_{i+j}, \ \  2 \leq i+j \leq
n-1 \\[2mm] e_n\circ e_1= \alpha e_{n-1},\quad e_n\circ e_n=\beta
e_{n-1} \\ \end{array}
\end{equation}
where  $\alpha, \beta \in \mathbb{C}$.
\end{pr}
\begin{proof} By theorem \ref{t2} we have that any $n$-dimensional
complex filiform Zinbiel algebra is isomorphic to an algebra of
the form:
$$
F_n+\beta
$$
where
$$ \beta(e_1,e_i)= 0, \ \mbox{ for } \
1 \leq i \leq n-1,
$$ $$
\beta(e_n,e_n) \in  lin\{e_3, e_4, \ldots,e_{n-1}\},
$$ $$
\beta(e_i,e_n), \  \beta(e_n,e_i) \in  lin\{e_{i+2}, e_{i+3},
\ldots, e_{n-1}\} \  \mbox{ for } \ 1 \leq i \leq n-1,
$$ $$
\beta(e_i,e_1) \in lin\{e_{i+2}, e_{i+3}, \ldots, e_{n-1}\} \
\mbox{ for } \  2\leq i \leq n-1,
$$ $$
\beta(e_i,e_j) \in lin\{e_{i+j+1}, e_{i+j+2}, \ldots ,e_{n-1}\} \
\mbox{ for } \ 2 \leq i \leq n-1.
$$

Similarly as in proof of the theorem \ref{t2} it is not difficult
to establish that the multiplications:
$$
e_i\circ e_1, \ \ e_i\circ e_j, \quad 1 \leq i, j \leq n-1
$$
can be obtained from $e_1\circ e_i=e_{i+1},$ $1 \leq i \leq n-2$
and identity (\ref{e1}).

By the similar procedure, we obtain
$$
\beta(e_i,e_1) =0,\  \mbox{ for } \  1 \leq i \leq n-1,
$$ $$
\beta(e_i,e_j)=0, \  \mbox{ for } \ 2 \leq i, j \leq
n-1.
$$

Let a Zinbiel algebra $A$ be isomorphic to the algebra
$F_n+\beta$. Then we have
$$
e_1\circ e_i=e_{i+1}, \  \mbox{ for } \  1\leq i \leq n-2.
$$

Put $e_1\circ e_n= \alpha_3e_3+ \alpha_4e_4+\ldots +
\alpha_{n-1}e_{n-1}$, $ e_n\circ e_1=\beta_3e_3+\beta_4e_4+
\ldots+ \beta_{n-1}e_{n-1}$.

Taking the change $e'_n=e_n- \alpha_3e_2-\alpha_4e_3- \ldots -
\alpha_{n-1}e_{n-2}$, we can suppose $\alpha_3= \alpha_4= \ldots =
\alpha_{n-1}=0$, i.e. $e_1 \circ e_n=0$.

Identity (\ref{e1}) implies
$$
(e_1\circ e_n)\circ e_1=e_1\circ(e_n\circ e_1)+ e_1\circ (e_1\circ
e_n) \Rightarrow e_1\circ (\beta_3e_3+\beta_4e_4+ \ldots
+\beta_{n-1}e_{n-1})=0
$$
therefore $\beta_3e_4+\beta_4e_5+ \ldots +\beta_{n-2}e_{n-1} =0$.
Hence, $e_n\circ e_1=\beta_{n-1}e_{n-1}$.

Consider the product
$$
(e_n\circ e_n)\circ e_1= e_n\circ (e_n \circ e_1)+ e_n\circ (e_1
\circ e_n)=0
$$
from which we obtain $e_n\circ e_n=\gamma e_{n-1}$ for some
$\gamma$.

Similarly to the proof of theorem \ref{t2} we can obtain that
$e_n\circ e_i=0$ for $2 \leq i \leq n-1$ and $e_i \circ e_n=0$ for
$1 \leq i \leq n-1$.

Thus, we obtain that Zinbiel algebra $A$ is split, i.e.
$A=NF_{n-1}\oplus C$. Applying theorem \ref{t1} complete the proof
of the proposition.
\end{proof}

The classification of complex filiform Zinbiel algebras is given
in the following

\begin{thm}\label{t3} An arbitrary $n$-dimensional ($n\geq 5$) complex filiform
Zinbiel algebras is isomorphic to one of the following pairwise
non isomorphic algebras:
$$
F_n^1: e_i\circ e_j=C_{i+j-1}^j  e_{i+j}, \quad 2 \leq i+j \leq
 n-1;
$$ $$
F_n^2: e_i\circ e_j=C_{i+j-1}^j  e_{i+j}, \quad 2 \leq i+j \leq
 n-1, \quad e_n\circ e_1=e_{n-1};
$$ $$
F_n^3: e_i\circ e_j=C_{i+j-1}^j  e_{i+j}, \quad 2 \leq i+j \leq
 n-1, \quad e_n\circ e_n=e_{n-1}.
$$
\end{thm}
\begin{proof} By proposition \ref{p1} we have the multiplication
(multiplication (\ref{e6})) in an $n-$dimensional complex filiform
Zinbiel algebra, namely,
$$
e_i\circ e_j=C_{i+j-1}^j  e_{i+j}, \quad 2 \leq i+j \leq n-1
$$ $$
e_n\circ e_1= \alpha e_{n-1}, \ \ e_n\circ e_n=\beta
e_{n-1}.
$$

Consider the general change of generator basic elements in the
form:
$$
e'_1=a_1e_1 + a_2e_2 + \ldots + a_ne_n ,
$$ $$
e'_n =b_1e_1 + b_2e_2 + \ldots + b_ne_n,
$$
where $a_1\neq o$ and $a_1b_n-a_nb_1\neq 0.$

Then expressing the rest basic elements of the new basis by means
of basis elements of the old basis and comparing the obtained
relations, we obtain the following restrictions:
$$a_1b_1=0,$$
$$a_1b_2+2a_2b_1=0,$$
$$ a_1b_3+ C_3^2 a_2b_2+3a_3b_1=0,$$
$$\dots \dots \dots \dots, $$
$$a_1b_{n-3}+ C_{n-3}^{n-4}a_2b_{n-4}+ C_{n-3}^{n-5}a_3b_{n-5}+ \ldots
+(n-3)a_{n-3}b_1=0,$$
$$ a_1b_{n-2}+ C_{n-2}^{n-3}a_2b_{n-3}+C_{n-2}^{n-4} a_3b_{n-4}+ \ldots
+(n-2)a_{n-2}b_1+\alpha a_nb_1+ \beta a_nb_n=0.$$

We have recurrently from these restrictions $b_1=b_2= \ldots
=b_{n-3}=0$. Hence, $b_{n-2}=-\frac{a_n b_n} {a_1} \beta$.

Consider the product
$$
e'_n \circ e'_1 =( -\frac{a_n b_n} {a_1} \beta
e_{n-2}+b_{n-1}e_{n-1}+b_ne_n) \circ (a_1e_1 + a_2e_2+ \ldots +
a_ne_n )=
$$
$$
= - (n-2)\beta a_n b_n e_{n-1}+ \alpha a_1 b_n e_{n-1}+\beta a_n
b_n e_{n-1}=(\alpha a_1b_n-(n-3) a_n b_n)e_{n-1}.
$$

On the other hand,
$$
e'_n\circ e'_1 = \alpha'a_1^{n-1}e_{n-1}.
$$

Comparing coefficients at $e_{n-1}$ we obtain
$$
\alpha a_1b_n-(n-3)\beta a_nb_n= \alpha'a_1^{n-1}.
$$

Consider the product

$ e'_n\circ e'_n =( -\frac{a_n b_n} {a_1} \beta
e_{n-2}+b_{n-1}e_{n-1}+b_ne_n)\circ ( -\frac{a_n b_n} {a_1} \beta
e_{n-2}+b_{n-1}e_{n-1}+b_ne_n)=b_n^2 \beta e_{n-1}.$

On the other hand,
$$
e'_n\circ e'_n= \beta' a_1^{n-1}e_{n-1}.
$$

Comparing coefficients we obtain
$$b_n^2\beta = \beta' a_1^{n-1}.$$

Now consider the following cases.

\textbf{ Case 1.} Let  $\beta=0$. Then  $\beta'=0$ and $\alpha b_n
= \alpha'a_1^{n-2}$. If $\alpha =0$, then $\alpha'=0$ and we have
algebra $F_n^1$. If $\alpha \neq 0$, then taking
$b_n=\frac{a_1^{n-2}} {\alpha}$ we get $\alpha'=1$, i.e. the
algebra $F_n^2$ is obtained.

\textbf{Case 2.} Let $\beta \neq 0$. Then putting
$b_n=\sqrt{a_1^{n-1}}$, $a_n=\frac{\alpha a_1} {(n-3)\beta}$ we
get $\beta'=1$, $\alpha'=0$ and algebra $F_n^3$ is obtained.

Note that the obtained algebras are not pairwise isomorphic.
\end{proof}

Comparing the description of complex filiform Leibniz algebras
\cite{Ayu} and the result of theorem \ref{t3} we can note how much
the class of complex filiform Zinbiel algebras is "thin". So,
although Zinbiel algebras and Leibniz algebras are Koszul dual,
but they quantitatively strongly are distinguished even in the
class of filiform algebras .

\section{Classification of four dimensional complex Zinbiel algebras.}

Since an arbitrary finite dimensional complex Zinbiel algebra is
nilpotent, therefore for an arbitrary four dimensional Zinbiel
algebra $A$ the condition $A^5=0$ holds.

As the direct sum of nilpotent Zinbiel algebras is nilpotent, then
in the description it is enough to consider non split algebras.

\begin{thm} An arbitrary four-dimensional complex non split Zinbiel
algebra is isomorphic to the one of the following pairwise non
isomorphic algebras:
$$
\begin{array}{ll}
A_1:  & e_1\circ e_1=e_2,  \ e_1\circ e_2=e_3,\  e_2\circ
e_1=2e_3,\  e_1\circ e_3=e_4, \  e_2\circ e_2=3e_4,  \\& e_3\circ
e_1=3e_4;\\ A_2: &  e_1\circ e_1=e_3,\  e_1\circ e_2=e_4, e_1\circ
e_3=e_4, \ e_3\circ e_1=2e_4; \\ A_3: &  e_1\circ e_1=e_3,\
e_1\circ e_3=e_4, \ e_2\circ e_2=e_4, \  e_3\circ e_1=2e_4;
\\ A_4:&  e_1\circ e_2=e_3,\  e_1\circ e_3=e_4,\  e_2\circ e_1=
-e_3; \\ A_5: & e_1\circ e_2=e_3,  \ e_1\circ e_3=e_4, \ e_2\circ
e_1= -e_3,\  e_2\circ e_2= e_4;\\ A_6: & e_1\circ e_1=e_4,\
e_1\circ e_2=e_3,\  e_2\circ e_1= -e_3,\  e_2\circ e_2= -2e_3+e_4;
\\  A_7: &  e_1\circ e_2=e_3,\  e_2\circ e_1=e_4,\  e_2\circ e_2=
-e_3; \\ A_8(\alpha):\ &  e_1\circ e_1=e_3,\  e_1\circ e_2=e_4,\
e_2\circ e_1= -\alpha e_3, e_2\circ e_2= -e_4, \ \alpha \in C; \\
A_9(\alpha): &  e_1\circ e_1=e_4,\  e_1\circ e_2=\alpha e_4,\
e_2\circ e_1= -\alpha e_4,\  e_2\circ e_2=e_4, \\ &  e_3\circ
e_3=e_4,\   \alpha \in C; \\ A_{10}: &  e_1\circ e_2=e_4,\
e_1\circ e_3=e_4,\ e_2\circ e_1= -e_4,\  e_2\circ e_2=e_4,\
e_3\circ e_1=e_4;\\ A_{11}: & e_1\circ e_1=e_4,\  e_1\circ
e_2=e_4,\  e_2\circ e_1= -e_4,\ e_3\circ e_3=e_4;\\ A_{12}: &
e_1\circ e_2=e_3,\  e_2\circ e_1=e_4;\\ A_{13}:& e_1\circ
e_2=e_3,\  e_2\circ e_1= -e_3, e_2\circ e_2=e_4;\\ A_{14}: &
e_2\circ e_1=e_4,\  e_2\circ e_2=e_3;\\ A_{15}(\alpha): & e_1\circ
e_2=e_4,\  e_2\circ e_2=e_3,\ e_2\circ e_1=\frac{1+\alpha}
{1-\alpha}e_4, \  \alpha \in C\setminus\{1\};\\ A_{16}: & e_1\circ
e_2=e_4,\  e_2\circ e_1= -e_4,\  e_3\circ e_3=e_4.\\
\end{array}
$$
\end{thm}
\begin{proof} Note that the result of Proposition 3.1 of \cite{Alb} also
holds for Zinbiel algebras. Therefore we have the following
possible cases for $(dim A^2,  dim A^3, dim A^4)$:
$$
(3, 2, 1),\ (2, 1, 0),\ (2, 0, 0),\ (1, 0, 0),\ (0, 0, 0).
$$

It is obvious, an Zinbiel algebra with the condition  $(3, 2, 1)$
is nul-filiform. Using theorem \ref{t1} we obtain algebra $A_1$.

Consider an algebra with the restriction $(2, 1, 0)$ (this algebra
is filiform).

Let $\{e_1, e_2, e_3, e_4\}$ be a basis of algebra $A$ satisfying
the conditions $A^2=\{e_3, e_4\}$, $A^3=\{e_4\}$. Then we can
suppose that
$$e_1\circ e_1=\alpha_1e_3+\alpha_2e_4,\quad e_1\circ
e_2=\alpha_3e_3+\alpha_4e_4,$$
$$e_2\circ e_1= \alpha_5e_3+\alpha_6e_4,\quad e_2\circ
e_2=\alpha_7e_3+\alpha_8e_4,$$
$$e_1\circ e_3=e_4,\quad e_2\circ e_3=\alpha_9e_4,$$
where $(\alpha_1, \alpha_3, \alpha_5,\alpha_7) \neq (0, 0, 0, 0)$.

\textbf{Case 1.} Let $(\alpha_1, \alpha_7)\neq (0, 0)$. Then by
arguments analogous to arguments in the proofs of theorems
\ref{t2} and \ref{t3} we obtain algebras:
$$e_1\circ e_1=e_3,\quad e_1\circ e_3=e_4,\quad e_3\circ e_1=2e_4;$$
$$ e_1\circ e_1=e_3,\quad e_1\circ e_2=e_4,\quad  e_1\circ e_3=e_4,\quad e_3\circ e_1=2e_4;$$
$$e_1\circ e_1=e_3,\quad e_1\circ e_3=e_4,\quad e_2\circ e_2=e_4,\quad e_3\circ
e_1=2e_4.$$

Note that the algebra defined by multiplication
$$e_1\circ e_1=e_3,\quad e_1\circ e_3=e_4,\quad e_3\circ e_1=2e_4$$
is split. So, in this case we have the algebras $A_2-A_3$.

\textbf{Case 2.} Let $(\alpha_1, \alpha_7)=(0, 0)$. Then
$(\alpha_3, \alpha_5)\neq (0, 0)$. If $\alpha_3\neq -\alpha_5$,
then taking $e'_1=Ae_1+e_2$, where $A\neq -\alpha_9$, we have case
1. It remains to consider the case $\alpha_3= -\alpha_5$. Denote
$e'_3=\alpha_3e_3+\alpha_4e_4$. Then we can write
$$e_1\circ e_1=\alpha_2e_4,\quad e_1\circ e_2=e_3,\quad e_2\circ
e_1=\alpha_5e_3+\alpha_6e_4,$$
$$e_2\circ e_2=\alpha_8e_4,\quad e_1\circ e_3=e_4,\quad e_2\circ e_3= \alpha_9e_4.$$

Consider the products
$$
(e_1\circ e_2)\circ e_1=e_1\circ (e_2\circ e_1)+e_1\circ (e_1\circ
e_2) = e_1\circ (-e_3+\alpha_6e_4)+e_1\circ e_3 =0 \Rightarrow
e_3\circ e_1=0,
$$ $$
(e_1\circ e_2)\circ e_2 = 2e_1\circ (e_2\circ e_2)=0 \Rightarrow
e_3\circ e_2=0.
$$

If we replace basic elements as follows:
$$e'_1=e_1, \ \ e'_2=e_2- \alpha_9e_1,\ \  e'_3=e_3- \alpha_2 \alpha_9 e_4,\ \  e'_4=e_4,$$
we obtain  $\alpha_9=0$, i.e. the multiplication in the algebra
has the following form:
$$e_1\circ e_1=\alpha e_4,\quad e_1\circ e_2=e_3,$$
$$e_2\circ e_1=-e_3+ \beta e_4,\quad e_2\circ e_2= \gamma e_4,\quad e_1\circ
e_3=e_4$$ (omitted products are equal to zero).

Check the isomorphism inside this family.

Consider the general change of generator basic elements:
$$
e'_1=a_1e_1+ a_2 e_2+a_3e_3,
$$ $$
e'_2=b_1e_1+ b_2e_2+b_3e_3
$$
where $a_1b_2-a_2b_1\neq 0$.

Expressing rest basic elements of the new basis via basic elements
of the old basis and comparing the obtained relations we obtain
the following restrictions:
$$a_1^2\alpha + a_1a_2\beta +a_2^2\gamma   + a_1a_3= \alpha'a_1^2
b_2,$$
$$ a_1b_2\beta +2a_2b_2\gamma +a_1b_3= \beta'a_1^2 b_2,$$
$$b_2\gamma  = \gamma'a_1^2, \ \ b_1=0.$$

Consider the following cases

\textbf{Case 2.1.} Let $\gamma=0$. Then $\gamma'=0$ and
 $$a_1\alpha +a_2\beta +a_3=\alpha'a_1 b_2,$$
 $$b_2 + b_3= \beta'a_1b_2.$$

Taking $a_3=-a_1\alpha-a_2\beta$ and $b_3=-b_2$ we obtain
$\alpha'=\beta'=0,$ i.e. we have the algebra $a_4.$

\textbf{Case 2.2.} Let  $\gamma \neq 0$. Then putting
$b_2=\frac{a_1^2} {\gamma}$, $a_3=-\frac{a_1^2\alpha + a_1
a_2\beta +a_2^2\gamma} {a_1},$ and $b_3=\frac{a_1b_2\beta+2a_2b_2
\gamma}{a_1},$ we get $\gamma'=1,$ $\alpha'=\beta'=0$, i.e. the
algebra $a_5$ is obtained.

Note that algebras with the conditions $(2, 0, 0), (1, 0, 0), (0,
0, 0)$ are associative. Therefore we can use the classification of
four dimensional algebras of Leibniz \cite{Alb}, i.e. choose
algebras with the condition $A^3=0$. \end{proof}

\end{document}